\documentclass[12pt]{article}

\usepackage{amsmath}
\usepackage{amssymb}

\begin{document}

\newcommand{\B}{{\mathcal B}}
\newcommand{\Bt}{{\wt{\mathcal B}}}
\newcommand{\SB}{{\mathcal S}{\mathcal B}}
\newcommand{\End}{{\rm{End}\ts}}
\newcommand{\non}{\nonumber}
\newcommand{\wt}{\widetilde}
\newcommand{\wh}{\widehat}
\newcommand{\ot}{\otimes}
\newcommand{\la}{\lambda}
\newcommand{\al}{\alpha}
\newcommand{\be}{\beta}
\newcommand{\ga}{\gamma}
\newcommand{\mut}{\wt{\mu}}
\newcommand{\hra}{\hookrightarrow}
\newcommand{\tpr}{t^{\tss\prime}}
\newcommand{\ve}{\varepsilon}
\newcommand{\ts}{\,}
\newcommand{\tss}{\hspace{1pt}}
\newcommand{\U}{ {\rm U}}
\newcommand{\Y}{ {\rm Y}}
\newcommand{\SY}{ {\rm SY}}
\newcommand{\C}{\mathbb{C}\tss}
\newcommand{\Z}{\mathbb{Z}}
\newcommand{\ZZ}{{\rm Z}}
\newcommand{\gl}{\mathfrak{gl}}
\newcommand{\sll}{\mathfrak{sl}}
\newcommand{\agot}{\mathfrak{a}}
\newcommand{\qdet}{ {\rm qdet}\ts}
\newcommand{\sdet}{ {\rm sdet}\ts}
\newcommand{\gr}{ {\rm gr}}
\newcommand{\sgn}{ {\rm sgn}}
\newcommand{\Sym}{\mathfrak S}

\newcommand{\Proof}{\noindent{\it Proof.}\ \ }   
\renewcommand{\theequation}{\arabic{section}.\arabic{equation}} 

\newtheorem{thm}{Theorem}[section]
\newtheorem{prop}[thm]{Proposition}
\newtheorem{cor}[thm]{Corollary}
\newtheorem{defin}[thm]{Definition}
\newtheorem{example}[thm]{Example}
\newtheorem{lem}[thm]{Lemma}

\newcommand{\bth}{\begin{thm}}
\renewcommand{\eth}{\end{thm}}
\newcommand{\bpr}{\begin{prop}}
\newcommand{\epr}{\end{prop}}
\newcommand{\ble}{\begin{lem}}
\newcommand{\ele}{\end{lem}}
\newcommand{\bco}{\begin{cor}}
\newcommand{\eco}{\end{cor}}
\newcommand{\bde}{\begin{defin}}
\newcommand{\ede}{\end{defin}}
\newcommand{\bex}{\begin{example}}
\newcommand{\eex}{\end{example}}

\newcommand{\bal}{\begin{aligned}}
\newcommand{\eal}{\end{aligned}}
\newcommand{\beq}{\begin{equation}}
\newcommand{\ben}{\begin{equation*}}

\def\beql#1{\begin{equation}\label{#1}}

\newbox\squ  
\setbox\squ=\hbox{\vrule width.3pt
            \vbox{\hrule height.3pt width.4em\kern1ex\hrule height.3pt}%
            \vrule width.3pt}
\def\endproof{%
  \ifmmode\eqno\copy\squ\smallskip
  \else{\unskip\nobreak\hfil%
    \penalty50\hskip2em\hbox{}\nobreak\hfil\copy\squ
    \parfillskip=0pt \finalhyphendemerits=0\penalty-100\smallskip}
  \fi}

\title{\Large\bf  Representations of reflection algebras}
\author{{A. I. Molev\quad and\quad E. Ragoucy}}

\date{} 
\maketitle

\vspace{7 mm}

\begin{abstract}
We study a class of algebras $B(n,l)$ associated with 
integrable models with boundaries. 
These algebras can be identified with coideal
subalgebras in the Yangian for $gl(n)$. We construct
an analog of the quantum determinant and show that its coefficients
generate the center of $B(n,l)$.
We develop an analog of Drinfeld's highest weight theory
for these algebras and
give a complete description
of their finite-dimensional irreducible representations.

\vspace{7 mm}


\end{abstract}


\vspace{30 mm}

\noindent
School of Mathematics and Statistics\newline
University of Sydney,
NSW 2006, Australia\newline
alexm@maths.usyd.edu.au

\vspace{7 mm}

\noindent
LAPTH, Chemin de Bellevue, BP 110\newline
F-74941 Annecy-le-Vieux cedex, France\newline
ragoucy@lapp.in2p3.fr

\newpage

\section{Introduction}\label{sec:int}
\setcounter{equation}{0}

A central role in the theory of integrable models
in statistical mechanics is played by the {\it Yang--Baxter equation\/}
\beql{ybeq}
R_{12}(u-v)\ts R_{13}(u-w)\ts R_{23}(v-w)=R_{23}(v-w)\ts R_{13}(u-w)\ts R_{12}(u-v);
\end{equation}
see Baxter~\cite{b:es}.
Here $R(u)$ is a linear operator $R(u): V\ot V\to V\ot V$ on the tensor square
of a vector space depending on
the {\it spectral parameter\/} $u$. Both sides of the Yang--Baxter equation
are linear operators on the triple tensor product $V\ot V\ot V$ and 
the indices of $R(u)$ indicate the copies of $V$ where $R(u)$ acts;
e.g., $R_{12}(u)=R(u)\ot 1$. A simplest nontrivial solution of the equation
is provided by the {\it Yang\/} $R$-{\it matrix\/}
\beql{yrm}
R(u)=1-P\ts u^{-1},
\end{equation}
where $P$ is the permutation operator $P:\xi\ot\eta\mapsto \eta\ot \xi$
in the space $\C^n\ot\C^n$.	The Yang $R$-matrix emerges in the {\it XXX\/} or
({\it six vertex\/})  {\it model\/} \cite{b:es}. 
It
gives rise to an algebra with the defining relations
given by the {\it RTT relation\/}
\beql{RTT}
R(u-v)\ts T_1(u)\ts T_2 (v) = T_2(v)\ts T_1(u)\ts R(u-v)
\end{equation}
(we discuss the precise meaning of this relation below in Section~\ref{sec:def}).
The algebraic structures associated with the Yang--Baxter equation were
studied in the works of Faddeev's school in the 70-s and 80-s
in relation with the {\it quantum inverse
scattering method\/}; see e.g. Takhtajan--Faddeev~\cite{tf:qi},
Kulish--Sklyanin~\cite{ks:qs}.
In particular, a central element called 
the {\it quantum determinant\/} in the algebra
defined by \eqref{RTT}
was introduced by 
Izergin and Korepin~\cite{ik:lm}
in the case of two dimensions.
The basic ideas and formulas associated with the
quantum determinant for an arbitrary $n$
are given in the paper Kulish--Sklyanin~\cite{ks:qs}. 
Tarasov~\cite{t:sq, t:im} described
irreducible representations ({\it monodromy matrices\/})
in the case $n=2$.

In was independently
realized by Drinfeld and Jimbo around 1985 that	
the algebraic structures associated 
with the quantum inverse scattering method
are naturally described by the language of 
{\it Hopf algebras\/}. This marked the beginning of the theory of
{\it quantum groups\/} \cite{d:qg}
(a historic background of this theory
is given in the book by Chari and Pressley~\cite{cp:gq}).
In his paper \cite{d:ha} Drinfeld introduced a remarkable class 
of quantum groups called the {\it Yangians\/}.
For any simple Lie algebra $\agot$ the Yangian $\Y(\agot)$
is a canonical deformation
of the universal enveloping algebra $\U(\agot[x])$ 
for the polynomial current
Lie algebra $\agot[x]$.	The Hopf algebra defined by the $RTT$ relation
\eqref{RTT} with the Yang $R$-matrix \eqref{yrm}
is called the {\it Yangian for the general linear Lie algebra\/}
$\gl(n)$ and denoted $\Y(\gl(n))$ or $\Y(n)$. 

The significance of the Yangians was explained by 
Drinfeld \cite{d:ha} who showed that the rational solutions
of the Yang--Baxter equation \eqref{ybeq} are described by
the Yangian representations.
The irreducible finite-dimensional representations 
were classified in his subsequent paper
\cite{d:nr}. The results
turned out to be parallel to those for the semisimple Lie algebras.
Every irreducible finite-dimensional representation of $\Y(\agot)$
is a quotient of the corresponding universal highest weight module
where the components of the highest weight satisfy some dominance
conditions.	Explicit constructions of all such representations
for the Yangians $\Y(\sll(2))$ and $\Y(2)$ are given in Tarasov~\cite{t:sq, t:im}
and Chari--Pressely~\cite{cp:yr}.
However, apart from this case, the explicit structure 
of the Yangian representations
remains unknown even in the case of $\gl(n)$ (a description of
a class of {\it generic\/} and {\it tame\/}
representations is given in \cite{m:gt} and \cite{nt:ry}
via Gelfand--Tsetlin bases).

Sklyanin \cite{s:bc} introduced a class of algebras
associated with integrable models with boundaries
(we call them the {\it reflection algebras\/} in this paper).
His approach
was inspired by Cherednik's scattering theory \cite{c:fp}
for factorized particles on the half-line. Instead of the 
$RTT$ relation the algebras are defined by the
{\it reflection equation\/}
\beql{re}
R(u-v)\ts B_1(u)\ts R(u+v)\ts B_2 (v) = B_2 (v)\ts 
R(u+v)\ts B_1(u)\ts R(u-v).
\end{equation}
In  \cite{s:bc}
commutative subalgebras in the reflection algebras 
(in the case of two dimensions) were constructed
and the {\it algebraic Bethe ansatz\/} was described.
Moreover, an analog of the quantum determinant
for these algebras was introduced and some properties
of the highest weight representations were discussed.

Different versions of \eqref{re} were employed in
the works Reshetikhin--Semenov-Tian-Shansky~\cite{rs:ce},
Olshanski~\cite{o:ty} and
Noumi~\cite{n:ms}.
Similar classes of algebras were studied in Kulish--Sklyanin~\cite{ks:as},
Kuznetsov--J\o{rgensen}--Christiansen~\cite{kjc:nb},
Koornwinder and Kuznetsov~\cite{kk:gh, k:hf}.
Recently, algebras of this kind were discussed in the physics literature
in connection with the NLS model, they
describe the integrals of motion of the model; 
see Liguori--Mintchev--Zhao~\cite{lmz:be},
Mintchev--Ragoucy--Sorba~\cite{mrs:ss}.

In this paper we consider a family of reflection algebras $\B(n,l)$.
They are defined as associative algebras whose generators
satisfy two types of relations:
the reflection equation and the unitary condition; see \eqref{quater} and 
\eqref{unitary}
below. The unitary condition allows us to identify $\B(n,l)$
with a subalgebra in the $\gl(n)$-Yangian $\Y(n)$; see Theorem~\ref{thm:emb}.
This condition was not explicitly used in \cite{s:bc}, but it appears e.g. in
\cite{kk:gh}. If we omit it we get
a larger algebra $\Bt(n,l)$ such that $\B(n,l)$ is a quotient
of $\Bt(n,l)$ by an ideal generated by some central elements.
The consequence of this fact is that
the finite-dimensional irreducible representations of both
algebras $\B(n,l)$ and $\Bt(n,l)$ are essentially the same.

On the other hand, the subalgebra $\B(n,l)$
turns out to be a (left) coideal
in the Hopf algebra $\Y(n)$; see Proposition~\ref{prop:coideal}. 
This allows us to regard
the tensor product $L\ot V$ of an $\Y(n)$-module $L$ and a $\B(n,l)$-module $V$
as a $\B(n,l)$-module; cf. \cite[Proposition~2]{s:bc}.

We show that the center of $\B(n,l)$ is generated by the coefficients
of an analog of the quantum determinant which we call, following \cite{mno:yc},
the {\it Sklyanin determinant\/}. We derive a formula which expresses
the Sklyanin determinant in terms of the 
quantum determinant for the Yangian $\Y(n)$.

The aforementioned properties of the algebras $\B(n,l)$ exhibit much
analogy with the {\it twisted Yangians\/} introduced by Olshanski~\cite{o:ty};
see also \cite{mno:yc} for a detailed exposition. Moreover,
in two dimensions the algebras $\B(2,0)$ and $\B(2,1)$
turn out to be respectively isomorphic to the symplectic and orthogonal twisted Yangians;
see Section~\ref{subsec:n=2}.
This analogy also extends to the representation theory; cf. \cite{m:fd}.
We prove here that, as for the twisted Yangians, the Drinfeld
highest weight theory \cite{d:nr} is applicable to the algebras $\B(n,l)$.
Every finite-dimensional irreducible
representation of $\B(n,l)$ is highest weight, and given an irreducible
highest weight module, we produce necessary and 
sufficient conditions for it
to be finite-dimensional. These conditions are 
expressed in terms of the Drinfeld polynomials in a way similar
to \cite{d:nr} and \cite{m:fd} for the Yangians and twisted Yangians.
Note, however, an essential difference: here {\it all\/} Drinfeld
polynomials must satisfy a symmetry condition; see Theorem~\ref{thm:repgen}.
In particular, all of them have even degree.

In conclusion, we would like to emphasize the common feature
of the three classes of algebras, the Yangian, twisted Yangians
and reflection algebras: the defining relations in all the 
cases can be presented
in a special matrix from. 
This allows special algebraic techniques 
(the so-called $R$-{\it matrix formalism})
to be used to describe the algebraic structure and study 
representations of these algebras. On the other hand, 
the close relationship of these `quantum' algebras with
the matrix Lie algebras leads to applications in the 
classical representation theory; see e.g. \cite{m:ya} and references therein.
The Yangian symmetries
have been found in various areas of physics including
the theory of integrable models in statistical mechanics,
conformal field theory,	quantum gravity. We note a surprising
connection of the Yangian and twisted Yangians 
with the finite $\mathcal W$-algebras recently discovered in
\cite{rs:yf, rs:yr, r:ty}; see also~\cite{br:rp}.

\medskip

This work was done during the first author's visit to the
{\it Laboratoire d'Annecy-le-Vieux de Physique Th\'eorique\/}, Annecy, France.
He would like to thank the {\it Laboratoire\/} for the warm hospitality.

\section{Definitions and preliminaries}\label{sec:def}
\setcounter{equation}{0}

Recall first the definition of the $\gl(n)$-{\it Yangian\/} 
$\Y(n)$; see e.g. \cite{ks:qs}, \cite{d:nr}.
We follow the notation from \cite{mno:yc} where
a detailed account of the properties of $\Y(n)$ is given.
The {\it Yangian\/} $\Y(n)$ is the
complex associative algebra with the
generators $t_{ij}^{(1)},t_{ij}^{(2)},\dots$ where $1\leq i,j\leq n$,
and the defining relations
\beql{defrel}
[t_{ij}(u),t_{rs}(v)]=\frac{1}{u-v}\Big(t_{rj}(u)t_{is}(v)-t_{rj}(v)t_{is}(u)
\Big),
\end{equation}
where
\beql{series}
t_{ij}(u) = \delta_{ij} + t^{(1)}_{ij} u^{-1} + t^{(2)}_{ij}u^{-2} +
\cdots \in \Y(n)[[u^{-1}]]
\non
\end{equation}
and $u$ is a formal (commutative) variable. 
Introduce the matrix
\beql{tmatrix}
T(u)=\sum_{i,j=1}^n t_{ij}(u)\ot E_{ij}
\in \Y(n)[[u^{-1}]]\ot \End\C^n,
\end{equation}
where the $E_{ij}$ are the standard matrix units.
Then the relations~(\ref{defrel}) are equivalent to
the single $RTT$ relation \eqref{RTT}.
Here $T_1(u)$ and $T_2(u)$ are regarded as elements of
$\Y(n)[[u^{-1}]]\ot \End\C^n\ot \End\C^n$, the subindex of $T(u)$
indicates to which copy of $\End\C^n$ this matrix corresponds, and
\beql{R(u)}
R(u)=1-Pu^{-1},\qquad P=\sum_{i,j=1}^n E_{ij}\ot E_{ji}\in
(\End\C^n)^{\otimes 2}.
\non
\end{equation}

Now we introduce the reflection algebras.
Fix a decomposition of the parameter $n$ into the sum
of two nonnegative integers, $n=k+l$.
Denote by $G$ the diagonal $n\times n$-matrix
\beql{G}
G=\text{diag}(\ve_1,\dots,\ve_n)
\end{equation}
where 
$\ve_i=1$ for $1\leq i\leq k$, and
$\ve_i=-1$  for $k+1\leq i\leq n$.
The {\it reflection algebra\/} $\B(n,l)$ is a unital associative algebra
with the generators $b_{ij}^{(r)}$ where $r$ runs over positive integers
and $i$ and $j$ satisfy $1\leq i,j\leq n$. To write down the defining
relations introduce formal series
\beql{biju}
b_{ij}(u)=\sum_{r=0}^{\infty}b_{ij}^{(r)}u^{-r}, \qquad b_{ij}^{(0)}=\delta_{ij}\ts\ve_i
\end{equation}
and combine them into the matrix 
\beql{bmatrix}
B(u)=\sum_{i,j=1}^n b_{ij}(u)\ot E_{ij}
\in \B(n,l)[[u^{-1}]]\ot \End\C^n.
\non
\end{equation}
The defining relations
are given by the {\it reflection equation\/}
\cite{s:bc}
\beql{quater}
R(u-v)\ts B_1(u)\ts R(u+v)\ts B_2 (v) = B_2 (v)\ts R(u+v)\ts B_1(u)\ts R(u-v)
\end{equation}
together with the {\it unitary condition\/}
\beql{unitary}
B(u)B(-u)=1,
\end{equation}
where we have used the notation of \eqref{RTT}.
Rewriting \eqref{quater} and \eqref{unitary} in terms of the matrix elements
we obtain, respectively,
\beql{defrelb}
\bal[{}]
[b_{ij}(u),b_{rs}(v)]&=\frac{1}{u-v}\Big(b_{rj}(u)b_{is}(v)-b_{rj}(v)b_{is}(u)\Big)\\
{}&+\frac{1}{u+v}\Big(\delta_{rj}\sum_{a=1}^n b_{ia}(u)b_{as}(v)-
\delta_{is}\sum_{a=1}^n b_{ra}(v)b_{aj}(u)\Big)\\
{}&-\frac{1}{u^2-v^2}\ts\delta_{ij}\Big(\sum_{a=1}^n b_{ra}(u)b_{as}(v)-
\sum_{a=1}^n b_{ra}(v)b_{as}(u)\Big)
\eal
\end{equation}
and
\beql{unitme}
\sum_{a=1}^n b_{ia}(u)b_{aj}(-u)=\delta_{ij}.
\end{equation}
We shall also be using the algebra $\Bt(n,l)$ which is defined
in the same way as $\B(n,l)$ but with the unitary condition \eqref{unitary}
dropped. (This corresponds to Sklyanin's original definition
\cite{s:bc}). We use the same notation $b_{ij}^{(r)}$
for the generators of $\Bt(n,l)$. As we shall see in the following proposition,
the algebra $\B(n,l)$ is isomorphic to a quotient of $\Bt(n,l)$ by
an ideal whose generators are central elements; see also \cite{mrs:ss}.

\bpr\label{prop:bbt}
In the algebra $\Bt(n,l)$ the product $B(u)B(-u)$ is a scalar matrix
\beql{bbt}
B(u)B(-u)=f(u)\ts 1,
\end{equation}
where $f(u)$ is an even series in $u^{-1}$
whose all coefficients are central in $\Bt(n,l)$.
\epr

\Proof Multiply both sides of \eqref{defrelb} by $u^2-v^2$
and put $v=-u$. We obtain
\beql{vmu}
\bal
{}&2u\ts\Big(\delta_{rj}\sum_{a=1}^n b_{ia}(u)b_{as}(-u)-
\delta_{is}\sum_{a=1}^n b_{ra}(-u)b_{aj}(u)\Big)\\
{}&= \delta_{ij}\Big(\sum_{a=1}^n b_{ra}(u)b_{as}(-u)-
\sum_{a=1}^n b_{ra}(-u)b_{as}(u)\Big).
\eal
\end{equation}
Choosing appropriate indices $i,j,r,s$, it is easy to see
that $B(u)B(-u)=B(-u)B(u)$ and that this matrix is scalar.
Thus, \eqref{bbt} holds for an even series $f(u)$.
Now multiply both sides of
\eqref{quater} by $B_2(-v)$ from the right:
\beql{quatf}
R(u-v)B_1(u)R(u+v)f(v) = B_2 (v)R(u+v)B_1(u)R(u-v)B_2(-v).
\end{equation}
Applying \eqref{quater} to the right hand side we write it as
\beql{quarhs}
\bal
B_2 (v)R(u+v)B_1(u)R(u-v)B_2(-v)&=B_2 (v)B_2(-v)R(u-v)B_1(u)R(u+v)\\
{}&=f(v)R(u-v)B_1(u)R(u+v).
\eal
\end{equation}
This shows that $f(u)$ is central. 
\endproof

We would like to comment on the relevance of the choice of the
initial matrix $G$ in the expansion
\beql{matexp}
B(u)=G+\sum_{r=1}^{\infty}B^{(r)}u^{-r};
\end{equation}
see also \cite{mrs:ss}. We could take $G$ to be an arbitrary nondegenerate matrix. However,
as the proof of Proposition~\ref{prop:bbt} shows, the reflection equation
implies that $G^2$ is a scalar matrix. Since $G$ is nondegenerate,
the scalar is nonzero. On the other hand, as can be easily seen, 
for any constant $c$ and any
nondegenerate matrix $A$
the transformations
\beql{trans}
B(u)\mapsto c\ts B(u)\qquad\text{and}\qquad B(u)\mapsto A B(u)A^{-1}
\end{equation}
preserve the reflection equation. Accordingly, the matrix $G$
is then transformed as $G\mapsto c\ts G$ and $G\mapsto AGA^{-1}$.
Therefore, we may assume that $G^2=1$ and using an appropriate 
matrix $A$ we can bring $G$ to the form \eqref{G}.

Note also that both the reflection equation and unitary condition 
are preserved by
the change of sign $B(u)\mapsto -B(u)$. This implies that
the algebras $\B(n,l)$ and $\B(n,n-l)$ are isomorphic.
In what follows we assume that the parameter $l$ satisfies
$0\leq l\leq n/2$.

It is immediate from the definitions of the algebras $\B(n,l)$ and
$\Bt(n,l)$ that given any formal series $g(u)\in 1+u^{-1}\C[[u^{-1}]]$
the mapping
\beql{autog}
B(u)\mapsto g(u)\ts B(u)
\end{equation}
is an automorphism of the algebra $\Bt(n,l)$. If $g(u)$
satisfies $g(u)\ts g(-u)=1$ then \eqref{autog}
is an automorphism of $\B(n,l)$.

\section{Algebraic structure of $\B(n,l)$}\label{sec:as}
\setcounter{equation}{0}

Here we show that each $\B(n,l)$ can be identified with
a coideal subalgebra in the Yangian $\Y(n)$ and prove an analog
of the Poincar\'e--Birkhoff--Witt theorem for the algebra $\B(n,l)$.
Then we describe its center using 
an analog of the quantum determinant.

\subsection{Embedding  $\B(n,l)\hra \Y(n)$}

Denote by $T^{-1}(u)$ the inverse matrix for $T(u)$; see \eqref{tmatrix}.
It can be easily seen that
the matrix $T^{-1}(-u)$ satisfies the $RTT$ relation \eqref{RTT}.
This implies that the mapping $T(u)\to T^{-1}(-u)$ defines an algebra
automorphism of the Yangian $\Y(n)$; cf. \cite{mno:yc}.	The following
connection between $\B(n,l)$ and $\Y(n)$ was indicated in \cite{s:bc}.

\bth\label{thm:emb}
The mapping
\beql{emb}
\varphi:B(u)\mapsto T(u)\ts G\ts T^{-1}(-u)
\end{equation}
defines an embedding of the algebra 
$\B(n,l)$ into the Yangian $\Y(n)$.
\eth

\Proof First we verify that $\varphi$ is an algebra homomorphism.
Denote the matrix $T(u)\ts G\ts T^{-1}(-u)$ by $\wt {B}(u)$.
We obviously have $\wt {B}(u)\wt {B}(-u)=1$ and so \eqref{unitary}
is satisfied. Further, we have
\begin{multline}\label{quatt}
R(u-v)\wt{B}_1(u)R(u+v)\wt{B}_2 (v)\\
= R(u-v)T_1(u)\ts G_1\ts T_1^{-1}(-u) R(u+v) T_2(v)\ts G_2\ts T_2^{-1}(-v).
\end{multline}
We find from the $RTT$ relation \eqref{RTT} that
\beql{coter}
T_1^{-1}(-u) R(u+v) T_2(v)=T_2(v)  R(u+v) T_1^{-1}(-u).
\end{equation}
Therefore, since $G_i$ commutes with $T_j(u)$ for $i\ne j$ we bring
the expression \eqref{quatt} to the form
\beql{inquatt}
R(u-v)T_1(u)T_2(v)\ts G_1 R(u+v)\ts G_2\ts T_1^{-1}(-u)\ts T_2^{-1}(-v).
\end{equation}
One easily verifies that
\beql{quag}
R(u-v)G_1R(u+v)G_2=	G_2 R(u+v) G_1	R(u-v).
\end{equation}
Applying \eqref{RTT} and \eqref{quag} we write \eqref{inquatt} as
\beql{in2quatt}
T_2(v)T_1(u)G_2\ts  R(u+v)\ts G_1\ts T_2^{-1}(-v)\ts T_1^{-1}(-u) R(u-v),
\end{equation}
which coincides with $\wt{B}_2 (v)R(u+v)\wt{B}_1(u)	R(u-v)$
due to \eqref{coter}. Thus, $\wt{B}(u)$ satisfies 
the reflection equation \eqref{quater}.

We now show that $\varphi$ has trivial kernel. The Yangian $\Y(n)$
admits two natural filtrations; see \cite{mno:yc}. Here we use the one
defined by setting $\deg_1 t_{ij}^{(r)}=r$. Similarly, we define
the filtration on the algebra $\B(n,l)$ by $\deg_1 b_{ij}^{(r)}=r$.
Let us first verify that the homomorphism $\varphi$
is filtration-preserving. By definition, the matrix elements of $\wt{B}(u)$
are expressed as
\beql{matt}
\wt{b}_{ij}(u)=\sum_{a=1}^n \ve_a t_{ia}(u) \tpr_{aj}(-u),
\end{equation}
where the $\tpr_{ij}(u)$ denote the matrix elements of $T^{-1}(u)$.
Inverting the matrix $T(u)$ we come to the following expression for
$\tpr_{ij}(u)$:
\beql{invt}
\tpr_{ij}(u)=\delta_{ij}+\sum_{k=1}^{\infty}(-1)^k
\sum_{a_1,\dots,a_{k-1}=1}^n
t_{ia_1}^{\tss\circ}(u)	t_{a_1a_2}^{\tss\circ}(u)
\cdots t_{a_{k-1}j}^{\tss\circ}(u),
\end{equation}
where $t_{ij}^{\tss\circ}(u)=t_{ij}(u)-\delta_{ij}$.
Taking the coefficient at $u^{-r}$ for $r\geq 1$ we get
\beql{invtc}
t_{ij}^{\tss\prime\tss(r)}=\sum_{k=1}^r(-1)^{k}
\sum_{a_1,\dots,a_{k-1}=1}^n
\sum_{r_1+\cdots+r_k=r}
t_{ia_1}^{(r_1)} t_{a_1a_2}^{(r_2)}
\cdots t_{a_{k-1}j}^{(r_k)},
\end{equation}
with the last sum taken over positive integers $r_i$.	Therefore,
we find from \eqref{matt} that the degree of $\wt{b}_{ij}^{\ts(r)}$
does not exceed	$r$. Hence, $\varphi$ is filtration-preserving
and it defines a homomorphism of the corresponding graded
algebras
\beql{homgr}
\gr_1\ts \B(n,l)\to \gr_1\ts\Y(n).
\end{equation}
Let $\overline{t}_{ij}^{\ts(r)}$ denote the image of
$t_{ij}^{(r)}$ in the $r$-th component of $\gr_1\ts\Y(n)$.
The algebra $\gr_1\ts\Y(n)$ is obviously commutative, and
it was proved in \cite[Theorem~1.22]{mno:yc} that
the elements $\overline{t}_{ij}^{\ts(r)}$ are 
its algebraically independent
generators.

On the other hand, by the defining relations \eqref{defrelb}, the algebra
$\gr_1\ts \B(n,l)$ is also commutative. Denote by
$\overline{b}_{ij}^{\ts(r)}$ the image of
$b_{ij}^{(r)}$ in the $r$-th component of $\gr_1\ts\B(n,l)$.
We find from the unitary condition
\eqref{unitme} that the elements
\beql{gengr}
\begin{aligned}
\overline{b}_{ij}^{\ts(2p-1)}&,\qquad 1\leq i,j\leq k 
\qquad\text{or}\qquad k+1\leq i,j\leq n, \\
\overline{b}_{ij}^{\ts(2p)}&,\qquad 1\leq i\leq k<j\leq n 
\qquad\text{or}\qquad 1\leq j\leq k<i\leq n,
\end{aligned}
\end{equation}
with $p$ running over positive integers, generate the algebra $\gr_1\ts\B(n,l)$.
Now, by \eqref{matt}, the image of $\overline{b}_{ij}^{\ts(r)}$
under the homomorphism \eqref{homgr} has the form
\beql{imbgr}
\overline{b}_{ij}^{\ts(r)}\mapsto \big(\ts\ve_i\ts(-1)^{r-1}+\ve_j\ts\big)\ts 
\overline{t}_{ij}^{\ts(r)} + \big(\ \dots\ \big),
\end{equation}
where $(\ \dots\ )$ indicates a linear combination of monomials
in the elements $\overline{t}_{ij}^{\ts(s)}$ with $s<r$.
This implies that the elements	\eqref{gengr} 
are algebraically independent which completes the proof. \endproof

The proposition implies the following analog of the
Poincar\'e--Birkhoff--Witt theorem for the algebra $\B(n,l)$.

\bco\label{cor:pbw}
Given any total ordering on the elements
\beql{pbw}
\begin{aligned}
{b}_{ij}^{\ts(2p-1)}&,\qquad 1\leq i,j\leq k 
\qquad\text{or}\qquad k+1\leq i,j\leq n, \\
{b}_{ij}^{\ts(2p)}&,\qquad 1\leq i\leq k<j\leq n 
\qquad\text{or}\qquad 1\leq j\leq k<i\leq n,
\end{aligned}
\non
\end{equation}
with $p=1,2,\dots$, the ordered monomials in these elements
constitute a basis of the algebra $\B(n,l)$. \endproof
\eco

Due to Proposition~\ref{thm:emb}, we may regard $\B(n,l)$
as a subalgebra of $\Y(n)$ so that the generators $b_{ij}(u)$
are identified with the elements $\wt{b}_{ij}(u)$ given by \eqref{matt}.
Recall that the Yangian $\Y(n)$ is a Hopf algebra with the coproduct 
$
\Delta : \Y(n)\to
\Y(n)\ot\Y(n) 
$ 
defined by
\beql{copr}
\Delta (t_{ij}(u))=\sum_{a=1}^n t_{i\tss a}(u)\ot
t_{aj}(u).
\end{equation}

\bpr\label{prop:coideal}
The subalgebra $\B(n,l)$ is a left coideal in $\Y(n)$:
\beql{coideal}
\Delta \big(\B(n,l)\big)\subseteq \Y(n)\ot \B(n,l).
\end{equation}
\epr

\Proof It suffices to show that the images of the generators of
$\B(n,l)$ under the coproduct are contained in $\Y(n)\ot \B(n,l)$.
Since $\Delta$ is an algebra homomorphism, the images
of the matrix elements $\tpr_{ij}(u)$ of $T^{-1}(u)$ are given by
\beql{coinv}
\Delta (\tpr_{ij}(u))=\sum_{a=1}^n \tpr_{aj}(u)\ot
\tpr_{i\tss a}(u).
\end{equation}
Therefore, using \eqref{matt} we calculate that
\beql{coprb}
\Delta (b_{ij}(u))=\sum_{a,c=1}^n 
t_{i\tss a}(u)\tpr_{cj}(-u)\ot b_{ac}(u),
\end{equation}
completing the proof. 
\endproof

\subsection{Sklyanin determinant}

The {\it quantum determinant\/} $\qdet T(u)$
of the matrix $T(u)$ is a formal series in $u^{-1}$,
\beql{qdetd}
\qdet T(u)=1+d_1u^{-1}+d_2u^{-2}+\cdots, \qquad d_i\in 	\Y(n),
\end{equation}
defined by
\beql{qdet}
\qdet T(u)=\sum_{p\in \Sym_n} \sgn\tss p\cdot t_{p(1)1}(u)\cdots
t_{p(n)n}(u-n+1);
\end{equation}
see Izergin--Korepin~\cite{ik:lm},
Kulish--Sklyanin~\cite{ks:qs}.
The elements $d_1,d_2,\dots$ 
are algebraically independent generators of
the center of the algebra $\Y(n)$; see e.g. \cite{mno:yc}
for the proof. The quantum determinant
can be equivalently defined in a $R$-matrix form
\cite{ks:qs}, \cite{c:qg}; see also \cite{mno:yc}.
Consider the tensor product
\beql{tenprn}
\Y(n)[[u^{-1}]]\ot \End\C^n\ot\cdots\ot \End\C^n
\end{equation}
with $n$ copies of $\End\C^n$.
For complex parameters $u_1,\ldots,u_n$ set
\beql{train}
R(u_1,\dots ,u_n) =(R_{n-1,n})(R_{n-2,n}R_{n-2,n-1}) \cdots (R_{1n}
\cdots R_{12}),
\end{equation}
where  we abbreviate $R_{ij}=R_{ij}(u_i-u_j)$, and the subindices
enumerate the copies of  $\End\C^n$ in \eqref{tenprn}.
If we specialize to 
\beql{specu}
u_i=u-i+1,\qquad i=1,\dots,n
\end{equation}
then
$R(u_1,\dots ,u_n)$ becomes the one-dimensional
anti-symmetrization operator $A_n$
in the space $(\End\C^n)^{\ot n}$. The quantum determinant $\qdet T(u)$
is defined by the relation
\beql{fundet}
A_n\ts T_1(u) \cdots T_n(u-n+1) = 
A_n\ts\qdet T(u).
\end{equation}

Now, using \eqref{quater} and the relation
\beql{ybe}
R_{ij}R_{ir}R_{jr}=R_{jr}R_{ir}R_{ij},
\end{equation}
(the Yang--Baxter equation)
we derive the following relation for the $B$-matrices:
\beql{fundb}
\bal
R(u_1,\dots ,u_n)B_1(u_1)\wt{R}_{12}\cdots \wt{R}_{1n}
B_2(u_2)\wt{R}_{23}\cdots \wt{R}_{2n}B_{3}(u_{3})\cdots 
\wt{R}_{n-1,n} B_{n}(u_{n})=&\\
B_{n}(u_{n}) \wt{R}_{n-1,n}\cdots B_{3}(u_{3})\wt{R}_{2n}\cdots \wt{R}_{23}
B_2(u_2) \wt{R}_{1n}\cdots 	\wt{R}_{12}B_1(u_1)
R(u_1,\dots ,u_n)&,
\eal
\end{equation}
where $\wt{R}_{ij}=R_{ij}(u_i+u_j)$. When we specialize the parameters
$u_i$ as in \eqref{specu}, the element \eqref{fundb} will be equal
to the product of the anti-symmetrizer $A_n$ and a series in $u^{-1}$
with coefficients in $\B(n,l)$. We call this series
the {\it Sklyanin determinant\/} and denote it by $\sdet B(u)$.
That is, $\sdet B(u)$ is defined by
\beql{sdet}
A_n\ts B_1(u)\wt{R}_{12}\cdots \wt{R}_{1n}
B_2(u-1)\cdots 
\wt{R}_{n-1,n} B_{n}(u-n+1)=A_n\ts\sdet B(u).
\end{equation}
As follows from the definition, the constant term of  $\sdet B(u)$ is
$\det G=(-1)^l$, so
\beql{sdetdec}
\sdet B(u)=(-1)^l+c_1\ts u^{-1}+c_2\ts u^{-2}+\cdots,\qquad c_i\in\B(n,l).
\end{equation}

In the next theorem we regard $\B(n,l)$ as a subalgebra in $\Y(n)$; 
see Theorem~\ref{thm:emb}. 

\bth\label{thm:center}
We have the identity
\beql{sqdet}
\sdet B(u)=\theta(u)\ts \qdet T(u)\ts \big(\qdet T(-u+n-1)\big)^{-1},
\end{equation}
where
\beql{theta}
\theta(u)=(-1)^l\ts\prod_{i=1}^k(2u-2n+2i)\ts\prod_{i=1}^l(2u-2n+2i)
\ts\prod_{i=1}^n
\frac{1}{2u-2n+i+1}.
\end{equation}
In particular, all the coefficients of $\sdet B(u)$ are central
in $\B(n,l)$. Moreover, the odd coefficients $c_1,c_3,\dots$
are algebraically independent and 
generate the center of the algebra $\B(n,l)$.
\eth

\Proof 	Substitute $B(u)=T(u)\ts G\ts T^{-1}(-u)$ into \eqref{sdet}.
Applying relation \eqref{coter} repeatedly, we bring
the left hand side of \eqref{sdet} to the form
\begin{multline}\label{sdett}
A_n\ts T_1(u)\cdots T_n(u-n+1)\ts G_1\ts \wt{R}_{12}\cdots \wt{R}_{1n}
\ts G_2\ts \cdots 
\wt{R}_{n-1,n}\ts  G_{n}\\
{}\times T^{-1}_1(-u)\cdots T^{-1}_n(-u+n-1).
\end{multline}
Now use \eqref{fundet}, and then note that
\beql{anthet}
A_n\ts G_1\ts \wt{R}_{12}\cdots \wt{R}_{1n}
\ts G_2\ts \cdots 
\wt{R}_{n-1,n}\ts  G_{n}=A_n\ts \theta(u),
\end{equation}
for some scalar function $\theta(u)$. Indeed, this follows e.g. from
\eqref{fundb} where we specialize $u_i$ as in \eqref{specu}
and consider the trivial representation of $\B(n,l)$ such that
$B(u)\mapsto G$. Furthermore, we have
\beql{qdinv}
A_n\ts T^{-1}_1(-u)\cdots T^{-1}_n(-u+n-1)=A_n\ts \big(\qdet T(-u+n-1)\big)^{-1}.
\end{equation}
This follows from \eqref{fundet}, if we first multiply both sides
by $T^{-1}_n(u-n+1)\cdots T^{-1}_1(u)$ from the right, replace $u$ with
$-u+n-1$ and then conjugate the left hand side by the permutation
of the indices $1,\dots,n$ which sends $i$ to $n-i+1$.

To complete the proof of \eqref{sqdet} we need to calculate $\theta(u)$.
It suffices to find a diagonal matrix element of the operator on
the left hand side of \eqref{anthet} corresponding
to the vector $e_1\ot\cdots\ot e_n$, where the $e_i$ denote
the canonical basis of $\C^n$.
We have $(n-1)!\ts A_n=A_n\ts A'_{n-1}$, where $A'_{n-1}$ is the anti-symmetrizer
in the tensor product of the copies of $\End\C^n$ corresponding
to the indices $2,\dots,n$. Note that $A'_{n-1}$ commutes with $G_1$. Furthermore,
we have the	identity
\beql{apr}
A'_{n-1}\ts \wt{R}_{12}\cdots \wt{R}_{1n}=\wt{R}_{1n}\cdots \ts \wt{R}_{12}	A'_{n-1},
\end{equation}
which easily follows from \eqref{ybe}. This allows us to use induction
on $n$ to find the matrix element and to
perform the calculation of $\theta(u)$ which is now
straightforward.

Using  \eqref{sqdet} we can conclude that all the coefficients
of $\sdet B(u)$ belong to the center of the Yangian $\Y(n)$, and hence
to the center of its subalgebra $\B(n,l)$.

Furthermore, if we put
\beql{sdetm}
c(u)=\sdet B\big(u+\frac{n-1}{2}\big)\ts \theta\big(u+\frac{n-1}{2}\big)^{-1},
\end{equation}
then by \eqref{sqdet} we have the identity
\beql{ssdet}
c(u)=d(u)\ts d(-u)^{-1},\qquad d(u)=\qdet T\big(u+\frac{n-1}{2}\big).
\end{equation}
The coefficients of $d(u)$ are algebraically independent generators of
the center of $\Y(n)$. Since $c(u)\ts c(-u)=1$, repeating
the argument of the second part of the proof of Theorem~\ref{thm:emb}
for the case of $\B(1,0)$ we find that all the even coefficients of $c(u)$
can be expressed in terms of the odd ones, 
and the latter are algebraically independent.
Since
\beql{sdetc}
\sdet B(u)=\theta(u)\ts c\big(u-\frac{n-1}{2}\big),
\end{equation}
the same holds for the coefficients of $\sdet B(u)$.

Finally, let us show that the center of $\B(n,l)$ is generated by the
coefficients of $\sdet B(u)$. We use another filtration on $\B(n,l)$
defined by setting $\deg_2 b_{ij}^{(r)}=r-1$. We first verify
that the corresponding graded algebra $\gr_2\ts\B(n,l)$ is isomorphic to
the universal enveloping algebra for a twisted polynomial current Lie algebra.
Consider the involution $\sigma$ of the Lie algebra $\gl(n)$ given
by $\sigma:E_{ij}\mapsto \ve_i\ts\ve_j\ts E_{ij}$. 
Denote by $\agot_0$ and $\agot_1$ the eigenspaces of $\sigma$
corresponding to the eigenvalues $1$ and $-1$, respectively.
In particular, $\agot_0$ is a Lie subalgebra of $\gl(n)$ isomorphic to 
$\gl(k)\oplus \gl(l)$. Denote by $\gl(n)[x]^{\sigma}$
the Lie algebra	of polynomials in a variable $x$ of the form 
\beql{polc}
a_0+a_1\ts x+a_2\ts x^2+\cdots+a_m\ts x^m,\qquad a_{2i}\in \agot_0,\quad 
a_{2i-1}\in \agot_1.
\end{equation}
We claim that the following is an algebra isomorphism:
\beql{gr2}
\gr_2\ts\B(n,l)\simeq \U(\gl(n)[x]^{\sigma}).
\end{equation}
Indeed, denote by
$\overline{b}_{ij}^{\ts(r)}$ the image of
$b_{ij}^{(r)}$ in the $(r-1)$-th component of $\gr_2\ts\B(n,l)$.
Then by the unitary condition we have for $r\geq 1$
\beql{gr2uni}
(\ve_i+(-1)^r\ts \ve_j)\ts \overline{b}_{ij}^{\ts(r)}=0,
\end{equation}
while the reflection equation gives
\beql{gr2qua}
[\overline{b}_{ij}^{\ts(r)}, \overline{b}_{kl}^{\ts(s)}]
=\delta_{kj}\ts (\ve_i\ts(-1)^{r-1}\ts +\ve_j)\ts\overline{b}_{il}^{\ts(r+s-1)}
-\delta_{il}\ts (\ve_i+\ve_j\ts (-1)^{r-1})\ts\overline{b}_{kj}^{\ts(r+s-1)}.
\end{equation}
This shows that the mapping
\beql{gr2map}
\overline{b}_{ij}^{\ts(r)}\mapsto (\ve_i+(-1)^{r-1}\ts \ve_j)\ts E_{ij}\ts x^{r-1}
\end{equation}
defines an algebra homomorphism $\gr_2\ts\B(n,l)\to\U(\gl(n)[x]^{\sigma})$.
Corollary~\ref{cor:pbw} ensures that its kernel is trivial.

Further, it is easily deduced from the definition \eqref{sdet}
of the Sklyanin determinant that the images of $\overline c_{2m+1}$
under the isomorphism \eqref{gr2} are
given by
\beql{imagec}
\overline c_{2m+1}\mapsto (-1)^l\ts 2\ts(E_{11}+\cdots+E_{nn})\ts x^{2m},
\qquad m\geq 0.
\end{equation}
The theorem will be proved if we show that the center of
$\U(\gl(n)[x]^{\sigma})$ is generated by
the elements $(E_{11}+\cdots+E_{nn})\ts x^{2m}$ with $m\geq 0$.
This is equivalent to the claim that the center
of $\U(\sll(n)[x]^{\sigma})$ is trivial. Here $\sll(n)[x]^{\sigma}$
is the Lie algebra of polynomials of the form \eqref{polc}
where now $\agot_0=\sll(n)\cap \big(\gl(k)\oplus\gl(l)\big)$.
However, this follows from a slight modification of a more general
result: see \cite[Proposition~4.10]{mno:yc}. Namely, if
$\agot$ is any Lie algebra and $\sigma$ is its involution,
we define $\agot[x]^{\sigma}$ as the Lie algebra of polynomials
of type \eqref{polc}. Then an argument similar to \cite{mno:yc} proves that
if the center of $\agot$ is trivial, and the $\agot_0$-module $\agot_1$
has no nontrivial invariant elements then the center of $\U(\agot[x]^{\sigma})$
is trivial. \endproof

\noindent
{\it Remarks\/}. (1) The algebra $\B(n,l)$ can also be regarded as
a deformation of the universal enveloping algebra
$\U(\gl(n)[x]^{\sigma})$. To see this, 
introduce the deformation parameter	$h$ and rewrite
the defining relations for $\B(n,l)$ in terms
of the re-scaled generators $b_{ij}^{\prime(r)}=	b_{ij}^{(r)}\ts h^{r-1}$.
These define a family of algebras $\B(n,l)_h$. If $h\ne 0$ the algebra
$\B(n,l)_h$ is isomorphic to $\B(n,l)$ while for $h=0$ one obtains
the universal enveloping algebra $\U(\gl(n)[x]^{\sigma})$. 

(2)	 It would be interesting to find an explicit formula for $\sdet B(u)$
in terms of the generators $b_{ij}(u)$. 
\endproof

Recall that the Yangian $\Y(\sll(n))$ for the special linear Lie algebra
$\sll(n)$ can be defined as the subalgebra of
$\Y(n)$ which consists of the elements stable under
all automorphisms of the form $T(u)\mapsto h(u)\ts T(u)$, where
$h(u)$ is a series in $u^{-1}$ with the constant term 1; see \cite{mno:yc}.
Then one has the tensor product decomposition
\beql{tprde}
\Y(n)={\text Z}(n)\ot \Y(\sll(n)),
\end{equation}
where ${\text Z}(n)$ is the center of $\Y(n)$.
Define the {\it special reflection algebra\/} $\SB(n,l)$ by
\beql{sbain}
\SB(n,l)=\B(n,l)\cap \Y(\sll(n)).
\end{equation}
In other words, $\SB(n,l)$ consists of the elements of $\B(n,l)$
which are stable under all automorphisms \eqref{autog}.
It is implied by  \eqref{tprde} (cf. \cite[Proposition~4.14]{mno:yc}), 
that the following decomposition holds
\beql{tprsba}
\B(n,l)={\text Z}(n,l)\ot \SB(n,l),
\end{equation}
where ${\text Z}(n,l)$ is the center of $\B(n,l)$.

\section{Representations of $\B(n,l)$}\label{sec:rep}
\setcounter{equation}{0}

Here we show that the Drinfeld highest weight theory \cite{d:nr}
applies to the representations of the algebras $\B(n,l)$;
see also \cite{m:fd} for the case of twisted Yangians.
We then give a complete description of the finite-dimensional irreducible
representations of 	$\B(n,l)$.

\subsection{Highest weight representations}\label{subsec:hwrep}

Recall first Drinfeld's classification results for
representations of the Yangian $\Y(n)$ \cite{d:nr}; see also \cite{m:fd}.
A representation $L$ of the Yangian $\Y(n)$
is called {\it highest weight\/} if there exists a nonzero vector 
$\xi\in L$ such that $L$ is generated by $\xi$,
\beql{trian}
\begin{aligned}
t_{ij}(u)\ts\xi&=0 \qquad \text{for} \quad 1\leq i<j\leq n, \qquad \text{and}\\
t_{ii}(u)\ts\xi&=\la_i(u)\ts\xi \qquad \text{for} \quad 1\leq i\leq n,
\end{aligned}
\end{equation}
for some formal series $\la_i(u)\in 1+u^{-1}\C[[u^{-1}]]$.
The vector $\xi$ is called the {\it highest vector\/}
of $L$ and the set $\la(u)=(\la_1(u),\dots,\la_n(u))$
is the {\it highest weight\/} of $L$.

Let $\la(u)=(\la_1(u),\dots,\la_n(u))$
be an $n$-tuple of formal series. Then there exists a unique,
up to an isomorphism, irreducible highest weight module	$L(\la(u))$
with the highest weight $\la(u)$.
Any finite-dimensional irreducible representation of $\Y(n)$
is isomorphic to $L(\la(u))$ for some $\la(u)$.
The representation $L(\la(u))$ is finite-dimensional if and only if
there exist monic polynomials $Q_1(u),\dots,Q_{n-1}(u)$ in $u$ such that
\beql{fdcondy}
\frac{\lambda_i(u)}{\lambda_{i+1}(u)}=\frac{Q_i(u+1)}{Q_{i}(u)},
\qquad i=1,\dots,n-1.
\end{equation}
These relations are analogs of the dominance conditions in the representation theory
of semisimple Lie algebras. Similar relations involving
monic polynomials are given by Drinfeld~\cite{d:nr}	to describe
finite-dimensional irreducible representations of the Yangian $\Y(\agot)$
for any simple Lie algebra $\agot$. The polynomials $Q_i(u)$ are called
the {\it Drinfeld polynomials\/} (usually denoted by $P_i(u)$).

Let us now turn to the reflection algebras. 
A representation $V$ of the algebra $\B(n,l)$
is called {\it highest weight\/} if there exists a nonzero vector 
$\xi\in V$ such that $V$ is generated by $\xi$,
\beql{trianb}
\begin{aligned}
b_{ij}(u)\ts\xi&=0 \qquad \text{for} \quad 1\leq i<j\leq n, \qquad \text{and}\\
b_{ii}(u)\ts\xi&=\mu_i(u)\ts\xi \qquad \text{for} \quad 1\leq i\leq n,
\end{aligned}
\end{equation}
for some formal series $\mu_i(u)\in \ve_i+u^{-1}\C[[u^{-1}]]$.
The vector $\xi$ is called the {\it highest vector\/}
of $V$ and the set $\mu(u)=(\mu_1(u),\dots,\mu_n(u))$
is the {\it highest weight\/} of $V$.

\bth\label{thm:fdhw}
Every finite-dimensional irreducible representation $V$ of the algebra $\B(n,l)$
is a highest weight representation.	Moreover, $V$ contains a unique
(up to a constant factor) highest vector.
\eth

\Proof Our approach is quite standard; cf. \cite[Proposition~12.2.3]{cp:gq}, \cite{m:fd}. 
Introduce the subspace of $V$
\beql{vo}
V^{\tss0}=\{\eta\in V\ |\ b_{ij}(u)\ts\eta=0,\qquad 1\leq i<j\leq n\}.
\end{equation}
We show first that $V^{\tss0}$ is nonzero. The defining relations
\eqref{defrelb} give
\beql{bij1}
[b_{ij}^{(1)},b_{rs}(u)]= (\ve_i+\ve_j)
\big(\delta_{rj}\ts b_{is}(u)-\delta_{is}\ts b_{rj}(u)\big).
\end{equation}
This implies that $V^{\tss0}$ can be equivalently defined as
\beql{voi}
V^{\tss0}=\{\eta\in V\ |\ b_{i,i+1}(u)\ts\eta=0,\qquad i=1,\dots,n-1\}.
\end{equation}
The operators $\ve_i\ts b_{ii}^{(1)}$ are pairwise commuting and so they have
a common eigenvector	$\zeta\ne0$ in $V$.
If we suppose that $V^{\tss0}=0$ then there exists an infinite sequence
of nonzero vectors in $V$, 
\beql{seqv}
\zeta,\quad b_{i_1,i_1+1}^{(r_1)}\ts\zeta,\quad 
b_{i_2,i_2+1}^{(r_2)}\ts b_{i_1,i_1+1}^{(r_1)}\ts\zeta,\quad \cdots.
\end{equation}
By \eqref{bij1} they all are eigenvectors for
the operators $\ve_i\ts b_{ii}^{(1)}$, $i=1,\dots,n$ of different weights.
Thus, they are linearly independent which contradicts to the assumption
$\dim V<\infty$. So, $V^{\tss0}$ is nontrivial.

Next, we show that all the operators $b_{rr}(u)$ preserve the subspace $V^{\tss0}$.
We use a reverse induction on $r$. For $r=n$ we see from \eqref{defrelb}
that if $i<j<n$ then $b_{ij}(u)\ts b_{nn}(v)\equiv 0$, 
where the equivalence is modulo the left ideal in $\B(n,l)$ generated by
the coefficients of $b_{ij}(u)$ with $i<j$.	Similarly, if $i<n$ then
the assertion follows from the equivalence
\beql{binn}
b_{in}(u)\ts b_{nn}(v)\equiv \frac{1}{u+v}\ts b_{in}(u)\ts b_{nn}(v)
\end{equation}
which is immediate from \eqref{defrelb}.
Now let $r<n$. Note the following obvious consequence of \eqref{defrelb}:
if $i<m$ and $i\ne j$ then
\beql{bijj}
b_{ij}(u)\ts b_{jm}(v)\equiv \frac{1}{u+v}\sum_{a=m}^n b_{ia}(u)\ts b_{am}(v).
\end{equation}
In particular, the right hand side is independent of $j$ and so, if $i<m$
then for any indices $j,j'\ne i$ we have
\beql{bikk}
b_{ij}(u)\ts b_{jm}(v)\equiv b_{ij'}(u)\ts b_{j'm}(v).
\end{equation}
If $i+1<r$ then it is immediate from \eqref{defrelb} that
$b_{i,i+1}(u)\ts b_{rr}(v)\equiv 0$. Further, using  \eqref{bijj}
we derive from \eqref{defrelb} that 
\beql{bk-1k}
b_{r-1,r}(u)\ts b_{rr}(v)\equiv \frac{n-r+1}{u+v}\ts b_{r-1,r}(u)\ts b_{rr}(v),
\end{equation}
which gives $b_{r-1,r}(u)\ts b_{rr}(v)\equiv 0$. Next, by \eqref{defrelb},
\beql{bkk+1}
\bal
b_{r,r+1}(u)\ts b_{rr}(v)
&\equiv \frac{1}{u-v}
\big( b_{r,r+1}(u)\ts b_{rr}(v)-  b_{r,r+1}(v)\ts b_{rr}(u)	\big)\\
{}&- \frac{1}{u+v} \sum_{a=r+1}^nb_{ra}(v)\ts b_{a,r+1}(u).
\eal
\end{equation}
By \eqref{bikk}, the second sum here is $-\dfrac{n-r}{u+v}\ts b_{r,r+1}(v)\ts b_{r+1,r+1}(u)$,
which is equivalent to zero by the induction hypothesis. So we get
\beql{bkk+1o}
\frac{u-v-1}{u-v}\ts b_{r,r+1}(u)\ts b_{rr}(v)+ \frac{1}{u-v}\ts 
b_{r,r+1}(v)\ts b_{rr}(u)\equiv 0.
\end{equation}
Swapping $u$ and $v$ we obtain
\beql{bkk+1uv}
-\frac{1}{u-v}\ts b_{r,r+1}(u)\ts b_{rr}(v)+ \frac{u-v+1}{u-v}\ts 
b_{r,r+1}(v)\ts b_{rr}(u)\equiv 0.
\end{equation}
The system of linear equations \eqref{bkk+1o} and  \eqref{bkk+1uv}
has only zero solution which proves the assertion in this case.
Finally, for $i>r$ we get from  \eqref{defrelb}
\beql{bik+1}
b_{i,i+1}(u)\ts b_{rr}(v)
\equiv \frac{1}{u-v}
\big( b_{r,i+1}(u)\ts b_{ir}(v)-  b_{r,i+1}(v)\ts b_{ir}(u)	\big),
\end{equation}
and
\beql{bii+1}
\bal
b_{r,i+1}(u)\ts b_{ir}(v)
&\equiv \frac{1}{u-v}
\big( b_{i,i+1}(u)\ts b_{rr}(v)-  b_{i,i+1}(v)\ts b_{rr}(u)	\big)\\
{}&- \frac{1}{u+v} \sum_{a=i+1}^nb_{ia}(v)\ts b_{a,i+1}(u).
\eal
\end{equation}
By \eqref{bikk}, the second sum here is equivalent to
$-\dfrac{n-i}{u+v}\ts b_{i,i+1}(v)\ts b_{i+1,i+1}(u)$,
which is equivalent to zero by the induction hypothesis. 
Therefore, \eqref{bii+1} implies that $b_{r,i+1}(u)\ts b_{ir}(v)$ is symmetric
in $u$ and $v$ which shows that
the left hand side of \eqref{bik+1} 
is equivalent to $0$.

Our next step is to show that all the operators $b_{rr}(u)$ on $V^{\tss0}$
commute.  For any $1\leq r\leq n$, the defining relations \eqref{defrelb} give
\beql{bkk}
\big(1-\frac{1}{u+v}\big)\ts[ b_{rr}(u), b_{rr}(v)]
\equiv \frac{1}{u+v}\ts \al_r(u,v),
\end{equation}
where
\beql{Akuv}
\al_r(u,v)=\sum_{a=r+1}^n \big(b_{ra}(u)\ts b_{ar}(v)-  b_{ra}(v)\ts b_{ar}(u)\big). 
\end{equation}
Using \eqref{defrelb} again we obtain for $a\geq r+1$
\begin{multline}\label{bkllk}
b_{ra}(u)\ts b_{ar}(v)-  b_{ra}(v)\ts b_{ar}(u)\\
\equiv \frac{1}{u+v}\ts\big([ b_{rr}(u), b_{rr}(v)]+[ b_{a\tss a}(u), b_{a\tss a}(v)]
+\al_r(u,v)+\al_r(u,v)\big).
\end{multline}
Taking the sum over $a$ gives
\beql{aksum}
(u+v-n+r)\al_r(u,v)\equiv (n-r)[ b_{rr}(u), b_{rr}(v)] +\sum_{a=r+1}^n
\big([ b_{a\tss a}(u), b_{a\tss a}(v)]+\al_r(u,v)\big).
\end{equation}
Using \eqref{bkk} we easily prove by a
reverse induction on $r$ that $[ b_{rr}(u), b_{rr}(v)]\equiv 0$ and
$\al_r(u,v)\equiv 0$. Finally, if $i<r$ then by \eqref{defrelb}
\beql{biikk}
[ b_{ii}(u), b_{rr}(v)]\equiv -
\frac{1}{u^2-v^2}\ts\big( [ b_{rr}(u), b_{rr}(v)]+\al_r(u,v)\big),
\end{equation}
which is equivalent to 0 as was shown above.

We can now conclude that the subspace $V^{\tss0}$ contains a 
common eigenvector $\xi\ne 0$ for the operators $b_{rr}(u)$,
that is, \eqref{trianb} holds for some formal series $\mu_i(u)$.

Since $V$ is irreducible, the submodule $\B(n,l)\ts\xi$ must coincide with $V$.
The uniqueness of $\xi$ now follows from Corollary~\ref{cor:pbw}. \endproof

Given any $n$-tuple $\mu(u)=(\mu_1(u),\dots,\mu_n(u))$, where
$\mu_i(u)\in \ve_i+u^{-1}\C[[u^{-1}]]$, we define
the {\it Verma module\/} $M(\mu(u))$ as the quotient of $\B(n,l)$ by
the left ideal generated by all the coefficients of the series $b_{ij}(u)$
with $1\leq i<j\leq n$, and $b_{ii}(u)-\mu_i(u)$ for $i=1,\dots,n$.
However, contrary to the case of the Yangian $\Y(n)$ the module
$M(\mu(u))$ can be trivial for some $\mu(u)$. If  $M(\mu(u))$ is nontrivial
we denote by $V(\mu(u))$ its unique irreducible quotient. Any
irreducible highest weight module with the highest weight $\mu(u)$
is clearly isomorphic to $V(\mu(u))$.

\bth\label{thm:exist}
The Verma module $M(\mu(u))$ is nontrivial (i.e. 
$V(\mu(u))$
exists) if and only if
\beql{mun}
\mu_n(u)\ts\mu_n(-u)=1,
\end{equation}
and for each $i=1,\dots,n-1$ the following conditions hold
\beql{muicond}
\wt{\mu}_i(u)\ts \wt{\mu}_i(-u+n-i)= \wt{\mu}_{i+1}(u)\ts \wt{\mu}_{i+1}(-u+n-i),
\end{equation}
where
\beql{mutil1}
\wt{\mu}_i(u)=(2u-n+i)\ts\mu_i(u)+\mu_{i+1}(u)+\cdots+\mu_n(u).
\end{equation}
\eth

\Proof	Suppose first that $V(\mu(u))$ exists. 	For each $i=1,\dots,n$ set
\beql{bbact}
\be_i(u,v)=\sum_{a=i}^n b_{ia}(u)\ts b_{ai}(v).
\end{equation}
The highest vector of $V(\mu(u))$ is an eigenvector of
$\be_n(u,v)$ with the eigenvalue $\mu_n(u)\ts\mu_n(v)$.
Due to the unitary condition \eqref{unitme} this eigenvalue must be equal to
1 if $v=-u$ which proves \eqref{mun}.
Using the notation of the proof of
Theorem~\ref{thm:fdhw},
we derive from \eqref{defrelb} that
\beql{bbdr}
\bal
\frac{u+v-n+i}{u+v}\ts \be_i(u,v)&\equiv b_{ii}(u)	
\ts b_{ii}(v)-\frac{1}{u+v}\sum_{a=i+1}^n
\be_a(v,u)\\
{}&+\frac{1}{u-v}\sum_{a=i+1}^n\big( b_{aa}(u)	
\ts b_{ii}(v)-b_{aa}(v)	\ts b_{ii}(u)\big).
\eal
\end{equation}
Therefore,
\beql{bedif}
\frac{u+v-n+i}{u+v}\big( \be_i(u,v)-\be_i(v,u)\big)
\equiv
\frac{1}{u+v}\sum_{a=i+1}^n\big( \be_a(u,v)-\be_a(v,u)\big).
\end{equation}
Since clearly $\be_n(u,v)\equiv \be_n(v,u)$, an easy induction implies that
$\be_i(u,v)\equiv \be_i(v,u)$ for all $i$.
Now consider \eqref{bbdr} with $i$ replaced by $i+1$ and subtract this from
\eqref{bbdr}. This gives
\beql{bdi}
\bal
\frac{u+v-n+i}{u+v}\ts \big(\be_i(u,v)-\be_{i+1}(u,v)\big)&\\
\equiv b_{ii}(u)	\ts b_{ii}(v)+
\frac{1}{u-v}&\sum_{a=i+1}^n\big( b_{aa}(u)	\ts b_{ii}(v)-b_{aa}(v)	\ts b_{ii}(u)\big)\\
{}- b_{i+1,i+1}(u)	\ts b_{i+1,i+1}(v)-
\frac{1}{u-v}&\sum_{a=i+2}^n\big( b_{aa}(u)	\ts b_{i+1,i+1}(v)-b_{aa}(v)	
\ts b_{i+1,i+1}(u)\big).
\eal
\end{equation}
Apply both sides to the highest vector of $V(\mu(u))$ and  put $u+v=n-i$.
We then get the following condition for the components of $\mu(u)$,
\beql{mudi}
\bal
\mu_i(u)	\ts \mu_{i}(v)+
\frac{1}{u-v}&\sum_{a=i+1}^n\big( \mu_{a}(u)	\ts \mu_{i}(v)-\mu_{a}(v)
\ts \mu_{i}(u)\big)\\
{}=\mu_{i+1}(u)	\ts \mu_{i+1}(v)+
\frac{1}{u-v}&\sum_{a=i+2}^n\big( \mu_{a}(u)	\ts \mu_{i+1}(v)-\mu_{a}(v)	
\ts \mu_{i+1}(u)\big),
\eal
\end{equation}
where $v=-u+n-i$. It is now a straightforward calculation to verify
that this condition is equivalent to \eqref{muicond}.

Conversely, suppose that the conditions \eqref{mun} and \eqref{muicond} hold.
We shall demonstrate that there exists a highest weight
module $L(\la(u))$ over the Yangian $\Y(n)$ such that $\B(n,l)$-cyclic span
of the highest vector is a $\B(n,l)$-module with the highest weight $\mu(u)$.
This will prove the existence of $V(\mu(u))$.

Suppose first that $L(\la(u))$ is an arbitrary 
irreducible highest weight module.
The {\it quantum comatrix\/}
$\widehat T(u)=(\wh t_{ij}(u))$ is defined by
\beql{comatrix}
\qdet T(u)=\widehat T(u)T(u-n+1).
\end{equation}
It can be deduced from \eqref{fundet} that
the matrix element $\widehat t_{ij}(u)$ 
equals $(-1)^{i+j}$ times the quantum determinant
of the submatrix of $T(u)$ obtained by removing the 
$i$th column and $j$th row. Therefore, the matrix element $\tpr_{ij}(u)$
of the inverse matrix $T^{-1}(u)$ can be expressed as
\beql{invme}
\tpr_{ij}(u)=	\big(\qdet T(u+n-1)\big)^{-1}\ts \widehat{t}_{ij}(u+n-1).
\end{equation}
This implies that the highest vector $\xi$ of the $\Y(n)$-module $L(\la(u))$
is annihilated by the elements $\tpr_{ij}(u)$ with $i<j$, and that
$\xi$ is an eigenvector for the $\tpr_{ii}(u)$. The corresponding eigenvalues
are found from the formula \eqref{qdet} and given by
\beql{eigenwt}
\tpr_{ii}(u)\ts\xi=\frac{\la_{i+1}(u+n-i)\cdots \la_n(u+1)}
{\la_{i}(u+n-i)\cdots \la_n(u)}\ts\xi.
\end{equation}
The following relations are easily derived
from \eqref{RTT}; see also \cite[Section~7]{mno:yc}
\beql{ttinv}
[t_{ij}(u),\tpr_{rs}(v)]=\frac{1}{u-v}
\Big(\delta_{rj}\sum_{a=1}^n t_{ia}(u)\ts\tpr_{as}(v)-
\delta_{is}\sum_{a=1}^n \tpr_{ra}(v)\ts t_{aj}(u)\Big).
\end{equation}
Hence modulo the left ideal in $\Y(n)$ 
generated by the elements $t_{ij}(u)$ with $i<j$ we can write for $a>i$
\beql{tti}
t_{ia}(u)\ts\tpr_{ai}(-u)\equiv \frac{1}{2u}\sum_{c=i}^n t_{ic}(u)\ts\tpr_{ci}(-u)
- \frac{1}{2u}\sum_{c=a}^n \tpr_{ac}(-u)\ts t_{ca}(u).
\end{equation}
Due to \eqref{matt} this implies that
\beql{biiy}
b_{ii}(u)\equiv \ve_i\ts t_{ii}(u)\ts\tpr_{ii}(-u)
+\frac{1}{2u}\sum_{a=i+1}^n\ve_a\Big(\sum_{c=i}^n t_{ic}(u)\ts\tpr_{ci}(-u)-
\sum_{c=a}^n \tpr_{ac}(-u)t_{ca}(u)\ts\Big).
\end{equation}
Now suppose that $i\geq k+1$ and set
\beql{fii}
f_{ii}(u)=-\sum_{c=i}^n \tpr_{ic}(-u)t_{ci}(u).
\end{equation}
Then \eqref{biiy} can be written as
\beql{biitr}
\frac{2u-n+i}{2u}\ts b_{ii}(u)\equiv	- t_{ii}(u)\ts\tpr_{ii}(-u)
-\frac{1}{2u}\sum_{a=i+1}^n f_{aa}(u).
\end{equation}
A similar transformation for $f_{ii}(u)$ gives
\beql{fiitr}
\frac{2u-n+i}{2u}\ts f_{ii}(u)\equiv	- \tpr_{ii}(-u)\ts	t_{ii}(u)
-\frac{1}{2u}\sum_{a=i+1}^n b_{aa}(u).
\end{equation}
Since $b_{nn}(u)\equiv f_{nn}(u)$, an easy induction proves that
$b_{ii}(u)\equiv f_{ii}(u)$ for all $i=k+1,\dots, n$.
By \eqref{biitr}, we have for those $i$, 
\beql{biitii}
\frac{2u-n+i}{2u}\ts b_{ii}(u)+\frac{1}{2u}\sum_{a=i+1}^n b_{aa}(u)
\equiv	- t_{ii}(u)\ts\tpr_{ii}(-u).
\end{equation}
In the same way for $i=1,\dots,k$ we get
\beql{biitlk}
\frac{2u-n+i}{2u-2l}\ts b_{ii}(u)+\frac{1}{2u-2l}\sum_{a=i+1}^n b_{aa}(u)
\equiv	t_{ii}(u)\ts\tpr_{ii}(-u).
\end{equation}
Now, applying both sides
of \eqref{biitii} and \eqref{biitlk} to $\xi$ and using \eqref{eigenwt}
we obtain
\beql{munla}
\mu_n(u)=
\begin{cases}\la_n(u)\la_n(-u)^{-1} \qquad&\text{if}\quad l=0,\\
-\la_n(u)\la_n(-u)^{-1} \qquad&\text{if}\quad l>0.
\end{cases}
\end{equation}
Moreover,
\beql{mupri}
\frac{\wt{\mu}_i(u)}{\wt{\mu}_{i+1}(u)}=
\frac{\la_i(u)\ts\la_{i+1}(-u+n-i)}{\la_{i+1}(u)\ts\la_i(-u+n-i)},
\end{equation}
if $i\ne k$; while
\beql{muprik}
\frac{\wt{\mu}_k(u)}{\wt{\mu}_{k+1}(u)}=\frac{l-u}{u}\cdot
\frac{\la_k(u)\ts\la_{k+1}(-u+l)}{\la_{k+1}(u)\ts\la_k(-u+l)}.
\end{equation}
The proof of the theorem is now completed as follows. By \eqref{mun} there exists
a series $\la_n(u)\in 1+u^{-1}\C[[u^{-1}]]$ such that \eqref{munla} holds.
Similarly, \eqref{muicond} ensures that there exist series
$\la_1(u),\dots,\la_{n-1}(u)$ satisfying
\eqref{mupri} and \eqref{muprik}. Thus, $V(\mu(u))$ is isomorphic to the
irreducible quotient of $\B(n,l)\ts\xi$. \endproof

\subsection{Representations of $\B(2,0)$ and $\B(2,1)$}\label{subsec:n=2}

The algebras $\B(2,0)$ and $\B(2,1)$ turn out to be isomorphic
to the symplectic and orthogonal {\it twisted Yangians\/} 
$\Y^-(2)$ and $\Y^+(2)$, respectively (see the Remark at the end of this section).
Thus the representations of $\B(2,0)$ and $\B(2,1)$ can be described
by using the results of \cite{m:fd} for the twisted Yangians.
In our argument we use the corresponding isomorphisms of the extended
algebras (Proposition~\ref{prop:isombtw}).
Note also that
this similarity between the reflection algebras and the twisted
Yangians does not extend to higher dimensions.

Recall (see e.g. \cite{mno:yc}) that the twisted Yangian $\Y^{\pm}(2)$ is 
an associative algebra with generators
$s_{ij}^{(1)},s_{ij}^{(2)},\dots$, where $i,j\in\{1,2\}$.
Introduce the generating series
\beql{defreltw}
s_{ij}(u)=\delta_{ij}+s_{ij}^{(1)}u^{-1}+s_{ij}^{(2)}u^{-2}+\cdots
\end{equation}
and combine them into the matrix $S(u)=(s_{ij}(u))$.
Then the defining relations have the form
of a reflection equation analogous to \eqref{quater}, 
\beql{quatw}
R(u-v)S_1(u)R^t(-u-v)S_2 (v) = S_2 (v)R^t(-u-v)S_1(u)R(u-v),
\end{equation}
as well as the symmetry relation
\beql{symmetry}
S^t(-u)=S(u)\pm\frac{S(u)-S(-u)}{2u}.
\end{equation}
Here we have used the notation of Section~\ref{sec:def}, and
\beql{Rtr}
R^t(u)=1-Qu^{-1},\qquad Q=\sum_{i,j=1}^2 E_{ij}^{\tss t}\ot E_{ji}\in
(\End\C^2)^{\otimes 2},
\non
\end{equation}
where the matrix transposition is defined by
\beql{transp}
A^t=\binom{a_{22}\quad a_{12}}{a_{21}\quad a_{11}}
\qquad\text{and}\qquad
A^t=\binom{a_{22}\quad -a_{12}}{-a_{21}\ \ \ \  a_{11}},
\end{equation}
for $\Y^{+}(2)$ and $\Y^{-}(2)$, respectively.

Another pair of associative algebras $\wt{\Y}^{+}(2)$ and $\wt{\Y}^{-}(2)$
is defined in the same way as $\Y^{\pm}(2)$, but with
the symmetry relation \eqref{symmetry} dropped. Then the algebra $\Y^{\pm}(2)$
is isomorphic to the quotient of $\wt{\Y}^{\pm}(2)$ by the ideal generated
by all coefficients of the series $\delta(u)$ defined by the relation
\beql{deltw}
\Big(1\mp\frac{1}{2u}\Big)\ts \delta(u)\ts Q= Q\ts S_1(u)\ts R(2u)\ts S_2^{-1}(-u).
\end{equation}
Moreover, these coefficients belong to the center of $\wt{\Y}^{\pm}(2)$;
see \cite[Theorems~6.3 and 6.4]{mno:yc} for the proofs of these assertions.
Note also that $\delta(u)$ satisfies $\delta(u)\ts\delta(-u)=1$
which is easy to deduce from \eqref{deltw}.

Recall that $\Bt(n,l)$ is the algebra with the defining relations
\eqref{quater}; see Section~\ref{sec:def}. Here we let $G={\rm diag\ts}(1,-1)$.

\bpr\label{prop:isombtw}
The mappings
\beql{epitwb}
S(u)\mapsto B\big(u+\frac12\big)\qquad\text{and}\qquad 	
S(u)\mapsto B\big(u+\frac12\big)\ts G
\end{equation}
define algebra isomorphisms $\wt{\Y}^{-}(2)\to \Bt(2,0)$ and
$\ \wt{\Y}^{+}(2)\to \Bt(2,1)$, respectively.
\epr

\Proof We have to show that the matrix $B\big(u+\frac12\big)$ or 
$B\big(u+\frac12\big)\ts G$,
respectively, satisfies the relation \eqref{quatw}.
Note that in the case of $\wt{\Y}^{-}(2)$ the operator
$P+Q$ is the identity on $(\End\C^2)^{\ot 2}$.
This implies that
\beql{reqrt}
R^t(-u-v)=\frac{u+v+1}{u+v}\ts R(u+v+1)
\end{equation}
proving the assertion. In the case of $\wt{\Y}^{+}(2)$ both $G_1\ts Q\ts G_1+P$
and $G_2\ts Q\ts G_2+P$ are the identity operators
which  implies that
\beql{reqrtp}
G_1\ts R^t(-u-v)\ts G_1=G_2\ts R^t(-u-v)\ts G_2=\frac{u+v+1}{u+v}\ts R(u+v+1).
\end{equation}
It remains to note that $G_1\ts G_2$ commutes with $R(u-v)$. \endproof  

Suppose now that $V$ is an irreducible finite-dimensional representation
of the twisted Yangian ${\Y}^{\pm}(2)$. Then $V$ is naturally extended
to $\wt{\Y}^{\pm}(2)$ and, by Proposition~\ref{prop:isombtw},
to the algebra $\Bt(2,0)$ or $\Bt(2,1)$, respectively.
The central element $f(u)$ defined in \eqref{bbt} acts in $V$
as a scalar series. Then there exists a series $g(u)\in 1+u^{-1}\C[[u^{-1}]]$
such that $g(u)\ts g(-u)\ts f(u)=1$ as an operator in $V$.
This means that the composition of $V$ with the corresponding
automorphism \eqref{autog} of $\Bt(2,0)$ or $\Bt(2,1)$ can be regarded
as an irreducible representation of the algebra $\B(2,0)$ or $\B(2,1)$, respectively.

Conversely, any irreducible finite-dimensional representation $V\!$ of
$\B(2,0)$ or $\B(2,1)$ is naturally extended to the corresponding algebra
$\Bt(2,0)$ or $\Bt(2,1)$ and hence to $\wt{\Y}^{\pm}(2)$.
The composition of $V$ with an appropriate automorphism $S(u)\mapsto h(u)\ts S(u)$
can be regarded as an irreducible representation of the twisted Yangian ${\Y}^{\pm}(2)$.
Indeed, the central element $\delta(u)$ defined by \eqref{deltw} 
acts as a scalar series on $V$,
and since $\delta(u)\ts\delta(-u)=1$,
the series $h(u)$ is found from the relation $h(u)\ts h(-u)^{-1}\ts \delta(u)=1$. 

This argument allows us to carry over
the description of representations of $\Y^{\pm}(2)$
to the case of the algebras $\B(2,0)$ or $\B(2,1)$. 
The following proposition is easily derived 
from \cite[Theorems~5.4 and 6.4]{m:fd}.
Suppose that two formal series $\mu_1(u)$ and $\mu_2(u)$ satisfy the conditions
of Theorem~\ref{thm:exist} so that the irreducible highest weight module
$V(\mu(u))$ exists.

\bpr\label{prop:repn=2}
{\em({\em i\ts})} The $\B(2,0)$-module $V(\mu(u))$ is finite-dimensional if and only if
there exists a monic polynomial $P(u)$ in $u$ such that $P(-u+2)=P(u)$
and
\beql{domb2}
\frac{(2u-1)\ts\mu_1(u)+\mu_2(u)}{2u\ts\mu_2(u)}=\frac{P(u+1)}{P(u)}.
\end{equation}
In this case $P(u)$ is unique.\par
{\em({\em ii\ts})} The $\B(2,1)$-module $V(\mu(u))$
is finite-dimensional if and only if there exist $\gamma\in\C$
and a monic polynomial $P(u)$ in $u$
such that $P(-u+2)=P(u)$, $P(\gamma)\ne 0$ and
\beql{domb1}
\frac{(2u-1)\ts\mu_1(u)+\mu_2(u)}{2u\ts\mu_2(u)}=
\frac{P(u+1)}{P(u)}\cdot\frac{\gamma-u}{\gamma+u-1}.
\end{equation}
In this case the pair $(P(u),\ga)$ is unique. \endproof
\epr

Given $\al,\be\in\C$ denote by $L(\al,\be)$ the irreducible 
$\gl(2)$-module with the highest weight $(\al,\be)$.
It is well-known that any finite-dimensional irreducible
representation of $\Y(2)$ is isomorphic to a tensor product of the form
\beql{tenpry2}
L=L(\al_1,\be_1)\ot\cdots\ot L(\al_k,\be_k),
\end{equation}
up to the twisting by an automorphism $T(u)\mapsto h(u)\ts T(u)$
of $\Y(2)$,
where $h(u)$ is a formal series \cite{t:im},
\cite{cp:yr}; see also \cite{m:fd}. 
The corresponding analogs of this result
for the algebras $\B(2,0)$ and $\B(2,1)$
are also immediate from \cite{m:fd} due to the above argument.
For any $\ga\in\C$ denote by $V(\ga)$ the one-dimensional representation
of $\B(2,1)$ such that the generators act by
\beql{1dim}
b_{11}(u)\mapsto\frac{u+\ga}{u-\ga},\qquad b_{22}(u)\mapsto -1,\qquad b_{12}(u)\mapsto 0,
\qquad b_{21}(u)\mapsto 0.
\end{equation}
Representations of this kind were constructed in the pioneering paper \cite{c:fp}.
Given any $\Y(n)$-module $L$ and any $\B(n,l)$-module $V$ we shall regard
$L\ot V$ as a $\B(n,l)$-module by using Proposition~\ref{prop:coideal}.

\bpr\label{prop:tenpr}
Any finite-dimensional irreducible 
$\B(2,0)$-module is isomorphic to the restriction of
a $\Y(2)$-module of the form \eqref{tenpry2}. 
Any finite-dimensional irreducible 
$\B(2,1)$-module is isomorphic to the tensor product
$L\ot V(\ga)$, where $L$ is a $\Y(2)$-module of the form \eqref{tenpry2}.
\endproof
\epr

\noindent
{\it Remarks\/.} (1) It is proved in \cite[Proposition~4.14]{mno:yc}
that the twisted Yangian admits the decomposition
\beql{dectwy}
\Y^{\pm}(2)=\ZZ^{\pm}\ot\SY^{\pm}(2),
\end{equation}
where $\ZZ^{\pm}$ is the center of $\Y^{\pm}(2)$, and $\SY^{\pm}(2)$
is a subalgebra called the {\it special twisted Yangian\/}.
Using the decomposition \eqref{tprsba} we conclude from
Proposition~\ref{prop:isombtw} that the algebras $\SY^{+}(2)$ and $\SY^{-}(2)$
are respectively isomorphic to $\SB(2,1)$ and $\SB(2,0)$. Indeed, each special
subalgebra is the quotient of the corresponding algebra $\wt{\Y}^{\pm}(2)$
or $\Bt(2,l)$ by the ideal generated by the center.	On the other hand,
both centers $\ZZ^{\pm}$ and $\ZZ(2,l)$ are polynomial algebras
in countably many variables; see Theorem~\ref{thm:center} above
and \cite[Theorem~4.11]{mno:yc}. This implies the isomorphisms
$\Y^{+}(2)\cong \B(2,1)$ and $\Y^{-}(2)\cong \B(2,0)$.

(2) Criteria of irreducibility of the modules $L$ and $L\ot V(\ga)$
(with $L$ given by \eqref{tenpry2}) 
over the algebras $\B(2,0)$ and $\B(2,1)$, respectively,
can also be obtained from the corresponding results in \cite{m:fd}.

\subsection{Classification theorem for $\B(n,l)$-modules}\label{subsec:class}

Theorem~\ref{thm:fdhw} together with the following result will
complete the description of finite-dimensional irreducible
representations of $\B(n,l)$.
Let $V(\mu(u))$ be an irreducible highest weight module over the algebra
$\B(n,l)$, so that the conditions of Theorem~\ref{thm:exist}
for the components $\mu_i(u)$ hold. We shall keep using
the notation \eqref{mutil1}.

\bth\label{thm:repgen}{\em({\em i\ts})}
The $\B(n,0)$-module $V(\mu(u))$ is finite-dimensional if and only if
there exist monic polynomials $P_1(u),\dots,P_{n-1}(u)$ in $u$ such that 
$P_i(-u+n-i+1)=P_i(u)$
and
\beql{dombgen0}
\frac{\mut_i(u)}{\mut_{i+1}(u)}=\frac{P_i(u+1)}{P_i(u)},\qquad i=1,\dots,n-1.
\end{equation}
\par
{\em({\em ii\ts})} The $\B(n,l)$-module $V(\mu(u))$ $($with $l>0$$\tss)$
is finite-dimensional if and only if there exist an element $\gamma\in\C$
and monic polynomials $P_1(u),\dots,P_{n-1}(u)$ in $u$ such that 
$P_i(-u+n-i+1)=P_i(u)$, $P_k(\gamma)\ne 0$ and
\begin{align}\label{dombgenl}
\frac{\mut_i(u)}{\mut_{i+1}(u)}&=\frac{P_i(u+1)}{P_i(u)},
\qquad i=1,\dots,n-1,\quad i\ne k,\\
\intertext{while}
\label{dombgenlk}
\frac{\mut_k(u)}{\mut_{k+1}(u)}&=\frac{P_k(u+1)}{P_k(u)}\cdot 
\frac{\gamma-u}{\gamma+u-l}. 
\end{align}
\eth

\Proof Set $V=V(\mu(u))$ and 
suppose first that $\dim V <\infty$. We shall use induction on
$n$. In the case $n=2$ the result holds by Proposition~\ref{prop:repn=2}.
Suppose now that $n\geq 3$ and consider the subspace
\beql{subspvt}
V_+=\{\eta\in V\ |\ b_{1i}(u)\ts\eta=0\qquad \text{for}\quad i=2,\dots,n\}.
\end{equation}
The calculations similar to those used in the proof of Theorem~\ref{thm:fdhw}
show that $V_+$ is stable under the operators $b_{ij}(u)$ with
$2\leq i,j\leq n$. Moreover, the operators 
\beql{repn-1}
b^{\circ}_{ij}(u)=b_{i+1,j+1}(u),\qquad i,j=1,\dots,n-1
\end{equation}
form a representation of the algebra $\B(n-1,l)$ in $V_+$ (recall that 
by our assumptions, $l\leq n/2$).
The cyclic span $\B(n-1,l)\ts \xi$ of the highest vector $\xi\in V$
is a module with the highest weight $\mu_+(u)=(\mu_2(u),\dots,\mu_n(u))$.
Since this module is finite-dimensional the conditions
of the theorem must be satisfied for the components of $\mu_+(u)$.

Similarly, consider the subspace
\beql{subspvtm}
V_-=\{\eta\in V\ |\ b_{in}(u)\ts\eta=0\qquad \text{for}\quad i=1,\dots,n-1\}.
\end{equation}
Now the operators $b_{nn}(u)$ and $b_{ij}(u)$ with $1\leq i,j\leq n-1$
preserve $V_-$. Moreover, the operators 
\beql{repn+1}
b^{\circ}_{ij}(u)=b_{ij}\big(u+\frac12\big)+\frac{\delta_{ij}}{2u}\cdot
b_{nn}\big(u+\frac12\big),
\qquad i,j=1,\dots,n-1
\end{equation}
form a representation of the algebra $\B(n-1,0)$ or $\B(n-1,l-1)$ in $V_-$,
respectively for $l=0$ and $l>0$. Again, the cyclic span
of $\xi$ is a finite-dimensional module over $\B(n-1,0)$ or $\B(n-1,l-1)$,
respectively, with the highest weight 
$\mu_-(u)=(\mu^{\circ}_1(u),\dots,\mu^{\circ}_{n-1}(u))$,
where 
\beql{muminco}
\mu^{\circ}_i(u)=\mu_i\big(u+\frac12\big)+\frac{1}{2u}\cdot\mu_n
\big(u+\frac12\big).
\end{equation}
By the induction hypothesis, the components of $\mu_-(u)$ must satisfy
the conditions of the theorem. Rewriting them in terms of the components of
$\mu(u)$ we complete the proof of the `only if\ts' part.

Conversely, suppose that the conditions of the theorem hold for the components of
the highest weight $\mu(u)$. Since $P_i(u)=P_i(-u+n-i+1)$ for each $i$,
there exist monic polynomials $Q_i(u)$ in $u$ such that
\beql{pqq}
P_i(u)=(-1)^{\deg Q_i}\ts Q_i(u)\ts Q_i(-u+n-i+1).
\end{equation}
Then there exists a representation $L(\la(u))$
of the Yangian $\Y(n)$ such that the components $\la_i(u)$ of $\la(u)$
satisfy the conditions
\eqref{fdcondy} for the polynomials $Q_i(u)$. The argument of the proof of
Theorem~\ref{thm:exist} shows that the $\B(n,0)$-span of the highest vector
$\xi$ of $L(\la(u))$ is a module with the highest weight $\mu(u)$
such that the relations \eqref{mupri} hold.
Since $V(\mu(u))$ is isomorphic to the irreducible quotient
of this span we conclude that $V(\mu(u))$ is finite-dimensional. 

Suppose now that $l>0$. It is easy to verify that for any $\ga\in\C$
the assignment
\beql{onedimn}
b_{ij}(u)\mapsto \delta_{ij}\ts \frac{u+\ga}{\ve_i\ts u-\ga}
\end{equation}
defines a one-dimensional representation of
$\B(n,l)$ which we denote by $V(\ga)$. 
Now consider the $\B(n,l)$-module $L(\la(u))\ot V(l-\ga)$; see Proposition~\ref{prop:coideal}.
Let $\eta$ be a basis vector of $V(l-\ga)$. 
Repeating again
the calculation of the proof of Theorem~\ref{thm:exist} we find that
the $\B(n,l)$-span of the vector $\xi\ot \eta$ is a module with the highest weight $\mu(u)$
such that the relations \eqref{mupri} hold for $i\ne k$, while
\beql{muprikga}
\frac{\wt{\mu}_k(u)}{\wt{\mu}_{k+1}(u)}=\frac{\ga-u}{\ga+u-l}\cdot
\frac{\la_k(u)\ts\la_{k+1}(-u+l)}{\la_{k+1}(u)\ts\la_k(-u+l)}.
\end{equation}
This implies that $V(\mu(u))$ is finite-dimensional.
\endproof

Using the decomposition \eqref{tprsba} we can deduce the following parametrization
of the representations of the special reflection algebra $\SB(n,l)$.

\bco\label{cor:sba}
{\em({\em i\ts})} The finite-dimensional irreducible representations of the
special reflection algebra $\SB(n,0)$ are in a one-to-one correspondence
with the families of monic polynomials $(P_1(u),\dots,P_{n-1}(u))$
such that $P_i(-u+n-i+1)=P_i(u)$.
\par
{\em({\em ii\ts})} The finite-dimensional irreducible representations of the
special reflection algebra $\SB(n,l)$ with $l>0$ are in a one-to-one correspondence
with the families $(P_1(u),\dots,P_{n-1}(u),\ga)$, where $\ga\in \C$
and the $P_i(u)$ are
monic polynomials such that
$P_i(-u+n-i+1)=P_i(u)$ and $P_k(\ga)\ne 0$.	 \endproof
\eco


\end{document}